\documentclass[11pt]{amsart}

\usepackage[latin1]{inputenc} 
\usepackage{amsmath,amsthm,amssymb,bm,mathtools,setspace} 
\usepackage{graphicx}
\usepackage{color}
\usepackage{verbatim}
\usepackage{float}
\usepackage{nccmath}
\usepackage[noadjust]{cite}
\usepackage[foot]{amsaddr}
\usepackage{enumitem}
\usepackage{thmtools}
\usepackage{thm-restate}
\usepackage{eqparbox}

\usepackage{pgfplots}
\pgfplotsset{width=10cm,compat=1.9}
\usepgfplotslibrary{fillbetween}

\setlist[enumerate,1]{label=\upshape(\roman*),ref=(\roman*)}

\newcommand{\COMMENT}[1]{}
\renewcommand{\COMMENT}[1]{\footnote{\textcolor{blue!70!black}{#1}}} 

\usepackage{hyperref}       
\usepackage[capitalise]{cleveref} 
\hypersetup{colorlinks=true,
   citecolor=blue,
   filecolor=blue,
   linkcolor=blue,
   urlcolor=blue
}

\creflabelformat{equation}{(#2#1#3)}
\crefdefaultlabelformat{#2#1#3}

\Crefname{enumi}{}{}
\Crefname{thm}{Theorem}{Theorems}
\Crefname{lm}{Lemma}{Lemmas}
\Crefname{cor}{Corollary}{Corollaries}
\Crefname{prop}{Proposition}{Propositions}
\Crefname{claim}{Claim}{Claims}
\Crefname{equation}{}{}
\Crefname{conjecture}{Conjecture}{Conjectures}
\Crefname{definition}{Definition}{Definitions}
\Crefname{figure}{Figure}{Figures}
\Crefname{fact}{Fact}{Facts}

\usepackage{array}  
\newcolumntype{C}{>{\centering \arraybackslash}m{3cm}}
\newcolumntype{D}{>{\centering \arraybackslash}m{1.25cm}}

\usepackage{enumitem} 
\setlist[enumerate]{topsep=0pt,itemsep=-0.5ex,partopsep=1ex,parsep=1ex}
\setlist[itemize]{topsep=0pt,itemsep=-0.5ex,partopsep=1ex,parsep=1ex}

\usepackage{bbm}

\theoremstyle{plain}
\newtheorem{theorem}{Theorem}[section]

\newtheorem{lemma}[theorem]{Lemma}
\newtheorem{cor}[theorem]{Corollary}

\theoremstyle{definition}

\theoremstyle{definition}

\renewcommand{\epsilon}{\varepsilon}

\renewcommand{\leq}{\leqslant}
\renewcommand{\geq}{\geqslant}

\usepackage[dvipsnames]{xcolor}

\usepackage{stackengine}
\stackMath
\newcommand{\stacklEq}[1]{%
  \setbox0=\hbox{${}\mathrel{\stackon[-1pt]{\leq}{\scriptstyle\text{#1\strut}}}{}$}
  \xdef\tmpwd{\dimexpr\the\wd0\relax}
  \kern.5\tmpwd\mathclap{\box0}&\kern.5\tmpwd
}

\newcommand{\stackEq}[1]{
  \setbox0=\hbox{${}\mathrel{\stackon[-1pt]{=}{\scriptstyle\text{#1\strut}}}{}$}
  \xdef\tmpwd{\dimexpr\the\wd0\relax}
  \kern.5\tmpwd\mathclap{\box0}&\kern.5\tmpwd
}

\allowdisplaybreaks

\newcommand*{\claimproofname}{Proof}

\newcommand{\ch}{\mathrm{ch}}
\newcommand{\fch}{\mathrm{Fch}}
\newcommand{\dem}{\mathrm{dem}}

\let\oldproofname=\proofname
\renewcommand{\proofname}{\rm\bf{\oldproofname}}

\title{A note on decompositions into 3-chromatic graphs}
\author[B. Granet]{Bertille Granet}
\email{first.lastname [at] warwick.ac.uk}
\address[Granet]{Warwick Mathematics Institute, University of Warwick, Coventry CV4 7AL, United Kingdom.}
\author[D. Horsley]{Daniel Horsley}
\email{first.lastname [at] monash.edu}
\address[Horsley, Santrio]{School of Mathematics, Monash University, Melbourne VIC 3800, Australia.}
\author[J. N. Santrio]{Jeffrey N. Santrio}
\email{first.lastname1 [at] monash.edu}
\thanks{Bertille Granet was partially supported by a Leverhulme Early Career Fellowship. Daniel Horsley was  supported by the Australian Research Council grant DP240101048. Jeffrey N. Santrio was supported by an Australian Government Research and Training Program Scholarship.}

\begin{document}
\setstretch{1.1}

\begin{abstract}
Delcourt and Postle recently introduced the `discharging in the dual' proof technique, leading to a breakthrough proof of Nash-Williams' conjecture on triangle decompositions. Independently, and almost simultaneously, the first and second authors built on the `condensation' method of Bryant, Dukes, Horsley, Maenhaut and Montgomery to determine the decomposition threshold of all odd cycles of length at least $73$.
In this note, we explain how the `discharging in the dual' and `condensation' methods can be combined to determine the decomposition threshold of cycles of any odd length, as well as the approximate decomposition threshold of a range of other 3-chromatic graphs.
\end{abstract}

\maketitle

\textbf{Note:} 
Shortly before our uploading this manuscript on arXiv, the AI-assisted preprint \cite{hou2026oddcycledecompositions} announced a proof of the 34 missing odd-cycle thresholds. Our work is independent of theirs. The key idea of scaling charges and then using \cite[Lemma~10.4]{DelPos2026} is the same in both pieces of work. However, our note present two main advantages. Firstly, we believe our proof is considerably clearer and more streamlined. Secondly, we thoroughly explore what the key method can be made to yield. Namely, we
fully investigate the range of weighted triangles and $3$-chromatic graphs to which it can be applied, and explain why it cannot work for the others.

\section{Introduction}

For a graph $F$, an \textit{$F$-decomposition} of a graph $G$ is a collection of copies of $F$ in $G$ whose edge sets partition the edge set of $G$. An obviously necessary condition for a graph $G$ to admit an $F$-decomposition is that $e(F)$ divides $e(G)$. Another one is that $\gcd\{d_F(u):u\in V(F)\}$ must divide $d_G(v)$ for every $v\in V(G)$. We say that $G$ is \textit{$F$-divisible} if the above two conditions are satisfied.

Given a graph $F$, we define the \textit{$F$-decomposition threshold}, denoted $\delta_F$, as the least real number $\delta\ge 0$ satisfying the following: for all $\delta'>\delta$, there exists $n_0\in \mathbb{N}$ such that whenever $G$ is an $F$-divisible graph on $n\ge n_0$ vertices with minimum degree $\delta(G)\ge \delta'n$, then $G$ has an $F$-decomposition.

The most notable result of our work is the following.

\begin{theorem}\label{T:main}
For each odd $\ell \geq 3$, we have $\delta_{C_\ell} = \frac{1}{2}+\frac{1}{2(\ell-1)}$. 
\end{theorem}

\cref{T:main} was previously proved for $\ell \geq 73$ in \cite{GraHor}, and for $\ell=3$ by Delcourt and Postle \cite{DelPos2026} in the course of proving the breakthrough result that every graph of sufficiently large order $n$ and minimum degree at least $\frac{3}{4}n$ has a $K_3$-decomposition if and only if it is $K_3$-divisible. We obtain \cref{T:main} from a more general result on approximate decomposition thresholds of 3-chromatic graphs, \cref{T:approx} below, by applying a powerful result of Barber, Kuhn, Lo and Osthus \cite{BarKuhLoOst}. To state \cref{T:approx}, we need to introduce some definitions.

A \emph{$k$-colouring} of a graph $F$ can be viewed as a partition of its vertex set into $k$ (possibly empty) independent sets called \emph{colour classes}. A graph is \emph{$k$-chromatic} if it has a $k$-colouring but no $(k-1)$-colouring. Let $\eta>0$. An \textit{$\eta$-approximate $F$-decomposition} of an $n$-vertex graph $G$ is an $F$-decomposition of a subgraph $G'$ of $G$ such that $|E(G')| \geq |E(G)|-\eta n^2$. We define $\delta_F^\eta$ as the least real number $\delta \ge 0$ satisfying the following: for all $\delta'>\delta$, there exists $n_0\in \mathbb{N}$ such that, if $G$ is a graph on $n\ge n_0$ vertices with $\delta(G)\ge \delta'n$, then $G$ has an $\eta$-approximate $F$-decomposition. The \textit{approximate $F$-decomposition threshold}, denoted $\delta_F^{0+}$, is defined as $\sup_{\eta >0}\delta_F^{\eta }$.

\begin{theorem}\label{T:approx}
Let $F$ be a $3$-chromatic graph, let $\{V_1,V_2,V_3\}$ be the colour classes of a $3$-colouring of $F$, and let $a \coloneqq e(V_1,V_2)$, $b \coloneqq e(V_1,V_3)$ and $c \coloneqq e(V_2,V_3)$ where $a \geq b \geq c$. Then
\begin{enumerate}
    \item
$\delta^{0+}_F \leq \frac{1}{2}+\max\{\frac{c}{2(a+b)},\frac{b+c}{8a}\}$.
    \item
$\delta^{0+}_F=\frac{1}{2}+\frac{c}{2(a+b)}$ if $c \geq \frac{b(a+b)}{3a-b}$ and no graph formed by removing fewer than $c$ edges from $F$ is bipartite.
\end{enumerate}
\end{theorem}

We note two interesting examples of families of graphs to which \cref{T:approx}(ii) can be applied. 

\begin{cor}\label{C:multi}
For $x,y,z \in \mathbb{Z}^+$ with $x \geq y \geq z \geq \frac{y(3y-z)}{y+z}$, we have $\delta^{0+}_{K_{x,y,z}}=\frac{1}{2}+\frac{yz}{2xy+2xz}$.
\end{cor}

\begin{cor}\label{C:wheel}
For odd $n \geq 5$, we have $\delta^{0+}_{W_n}=\frac{2}{3}$ where $W_n$ is the wheel with $n$ vertices.
\end{cor}

\cref{T:approx} is proved by combining the `condensation' methods introduced in \cite{BryDukHorMaeMon} with a result, \cref{T:weightedTriDecomp} below,
concerning fractional decompositions into weighted triangles. We define these concepts now. A \emph{weighted graph} $W$ is a pair consisting of a graph $F$, called the \textit{underlying graph}, and an edge-weighting $w:E(F)\to \mathbb{R}^+$. We let $V(W)\coloneq V(F)$ and $E(W)\coloneq E(F)$.
We define $T_{a,b,c}$ to be the weighted graph with underlying graph $K_3$ whose edges get assigned weight $a$, $b$, and $c$ respectively.

Let $W$ be a weighted graph, let $G$ be a graph, and let $I$ be the set of all injective homomorphisms from $W$ to $G$ (i.e.\ injective functions $i: V(W) \rightarrow V(G)$ such that $\{i(x),i(y)\} \in E(G)$ for all $\{x,y\} \in E(W)$). Then, a \emph{fractional $W$-decomposition of $G$} is an indexed set $\{s_i\}_{i \in I}$ of nonnegative real numbers such that, for each $e' \in E(G)$, we have $\sum_{e \in E(W)}\sum_{i \in I:i(e)=e'}s_iw(e)=1$. 

Given a weighted graph $W$, we define the \textit{fractional $W$-decomposition threshold}, denoted $\delta_W^*$, as the least real number $\delta\ge 0$ satisfying the following: for all $\delta'>\delta$, there exists $n_0\in \mathbb{N}$ such that whenever $G$ is a graph on $n\ge n_0$ vertices with $\delta(G)\ge \delta'n$, then $G$ has a fractional $W$-decomposition.

\begin{theorem}\label{T:weightedTriDecomp} 
Let $a,b,c \in \mathbb{R}^+$ such that $a \geq b \geq c$.
\begin{enumerate}
    \item
Any $n$-vertex graph $G$ with minimum degree at least $(\frac{1}{2}+\max\{\frac{c}{2(a+b)},\frac{b+c}{8a}\})n$ has a fractional $(T_{a,b,c})$-decomposition.
    \item
If $c \geq \frac{b(a+b)}{3a-b}$, then $\delta^*_{T_{a,b,c}}=\frac{1}{2}+\frac{c}{2(a+b)}$.
\end{enumerate}  
\end{theorem}

Upper bounds on $\delta_{T_{a,b,c}}^*$ were also established in \cite{BryDukHorMaeMon}. The bound of \cref{T:weightedTriDecomp} gives a strict improvement on these. We prove \cref{T:weightedTriDecomp} using the `discharging in the dual' method developed in \cite{DelPos2026}.

\section{Preliminaries}

Barber, K\"uhn, Lo, and Osthus \cite{BarKuhLoOst} used the method of `iterative absorption' to relate approximate decomposition thresholds to decomposition thresholds. The following result of theirs on odd cycles will allow us to derive \cref{T:main} from \cref{T:approx}.

\begin{theorem}[\cite{BarKuhLoOst}]\label{T:approxToIntegral}
For each odd $\ell \geq 3$, we have $\delta_{C_\ell}=\delta_{C_\ell}^{0+}$.
\end{theorem}

To derive \cref{T:approx} from \cref{T:weightedTriDecomp} we use the method of `condensations' introduced in \cite{BryDukHorMaeMon}. Let $F$ be a graph and let $\{V_1,\dots,V_k\}$ be the colour classes of a $k$-colouring of $F$. For distinct $i,j \in \{1,\ldots,k\}$, let $e_{ij}$ be the number of edges with one endpoint in $V_i$ and the other in $V_j$. The \textit{condensation of $F$ with respect to $V_1,\dots,V_k$} is the weighted graph $W$ with $V(W)\coloneq \{u_1,\dots,u_k\}$ such that, for all distinct $i,j \in \{1,\ldots,k\}$, the edge $u_iu_j$ is absent if $e_{ij}=0$ and is present with weight $e_{ij}$ otherwise. 
The key result on condensations we require is the following, which is a special case of \cite[Theorem 1.3]{BryDukHorMaeMon}.

\begin{theorem}[\cite{BryDukHorMaeMon}]\label{T:condensation}
If $F$ is a graph and $W$ is a condensation of $F$, then $\delta_{F}^{0+} \leq \delta_{W}^*$.
\end{theorem}

We introduce the concepts relating to charge assignments and duality that will be necessary for our proof of \cref{T:weightedTriDecomp}. A \emph{charge assignment} to a graph $G$ is a function $\phi: E(G) \rightarrow \mathbb{R}$. The \emph{total charge} of $\phi$ is $\sum_{e \in E(G)}\phi(e)$. Under $\phi$, an edge $e \in E(G)$ is \emph{positive}, \emph{negative} or \emph{zero} according to the sign of $\phi(e)$. We say that $\phi$ is \emph{$(T_{a,b,c})$-covered} if $a\phi(e_1)+b\phi(e_2)+c\phi(e_3) \geq 0$ for each ordered triple $(e_1,e_2,e_3)$ of three distinct edges of $G$ that form a triangle. Following \cite{DelPos2026}, we say that $\phi$ is \emph{triangle-covered} if it is $(T_{1,1,1})$-covered.

By definition, a fractional $(T_{a,b,c})$-decomposition of a graph $G$ exists if and only if there is a solution to a matrix equation of the form $A\mathbf{s}=\mathbf{1}$. Here, $A$ is a matrix whose rows correspond to the edges of $G$, whose columns correspond to the injective homomorphisms from $T_{a,b,c}$ to $G$, and in which the $(e,i)$-entry is the weight of the edge of $T_{a,b,c}$ that $i$ maps to $e$, or is 0 if no such edge exists. Applying Farkas' lemma \cite{Far} to this matrix equation gives the following.

\begin{lemma}\label{L:FarkasForConds}
A graph $G$ admits a fractional $(T_{a,b,c})$-decomposition if and only if each $(T_{a,b,c})$-covered charge assignment to $G$ has nonnegative total charge.
\end{lemma}

Lastly we need some lower bounds on fractional decomposition thresholds of weighted triangles. The following result follows from \cite[Lemma 3.1 and Lemma 3.3]{BryDukHorMaeMon}, but we give a short proof here for completeness and to highlight the corresponding constructions.

\begin{lemma}\label{C:LBforTris}
Let $a,b,c \in \mathbb{R}^+$ such that $a \geq b \geq c$. Then, $\delta^*_{T_{a,b,c}} \geq \max\{B(a,b,c),U(a,b,c)\}$, where
\[B(a,b,c)\coloneqq \mfrac{1}{2}+\mfrac{c}{2(a+b)} \quad\text{and}\quad U(a,b,c)\coloneqq \mfrac{1}{4}\left(3-\sqrt{\mfrac{2a-b-c}{2a+b+c}}\right).\]
\end{lemma}

Let $\delta \in (\frac{1}{2},\frac{3}{4}]$ be a rational number and let $n$ be a positive even integer such that $\delta n$ is also an integer.
Let $\mathcal{B}_n(\delta)$ be the family of $n$-vertex $\delta n$-regular graphs that contain a copy of $K_{\frac{n}{2},\frac{n}{2}}$; and let $\mathcal{U}_n(\delta)$ be the family of complete 4-partite graphs with two parts of size $(1-\delta)n$ and two parts of size $(\delta-\frac{1}{2})n$.
(We use $B$ and $\mathcal{B}$ because they correspond to graphs containing balanced complete bipartite graphs. Likewise $U$ and $\mathcal{U}$ correspond to unbalanced complete 4-partite graphs.)
These two families give rise to the lower bounds in \cref{C:LBforTris}.

\begin{proof}[Proof of \cref{C:LBforTris}]
We show that for any $\delta$ arbitrarily close to $\max\{B(a,b,c),U(a,b,c)\}$, there is an $n$-vertex graph with minimum degree at least $\delta n$ that does not have a fractional $(T_{a,b,c})$-decomposition. By \cref{L:FarkasForConds}, this can be done by showing the graph has a $(T_{a,b,c})$-covered charge assignment with negative total charge. 

Let $\delta \in (\frac{1}{2},\frac{3}{4}]$ be a rational number and let $n$ be a positive even integer such that $\delta n$ is also an integer. Graphs in $\mathcal{B}_n(\delta)$ are $(T_{a,b,c})$-covered if we assign charge $-\frac{c}{a+b}$ to each edge in the $K_{\frac{n}{2},\frac{n}{2}}$ and charge $1$ to each other edge. However, an easy calculation shows that this assignment will have negative total weight for any $\delta < B(a,b,c)$. This establishes the lower bound $\delta^*_{T_{a,b,c}} \geq B(a,b,c)$.

Graphs in $\mathcal{U}_n(\delta)$ are $(T_{a,b,c})$-covered if we assign charge $-\frac{b+c}{a}$ to each edge between two parts of the same size and charge $1$ to each other edge. However, a routine calculation shows that this assignment will have negative total weight for any $\delta < U(a,b,c)$. This establishes $\delta^*_{T_{a,b,c}} \geq U(a,b,c)$. 
\end{proof}

\section{Reductions of charge assignments}

In order to apply \cref{L:FarkasForConds}, we need a way to show that a $(T_{a,b,c})$-covered charge assignment has nonnegative total charge. To do this, we employ a slight generalisation of the discharging methods used by Delcourt and Postle in \cite{DelPos2026}. In order to avoid having to reprove many of their results in the general setting of $(T_{a,b,c})$-covered graphs (rather than just the triangle-covered setting), we employ a trick involving scaling the positive charges. This will allow us to use a key result of \cite{DelPos2026} directly.

Given $\alpha \in (0,1]$ and a charge assignment $\psi$ to a graph $G$, we define the \emph{$\alpha$-reduction} of $\psi$ to be the  charge assignment $\psi'$ to $G$ given by, for each $e \in E(G)$,
\[
\psi'(e) \coloneqq
\left\{
  \begin{array}{ll}
    \alpha\,\psi(e) & \hbox{if $\psi(e) \geq 0$} \\
    \psi(e) & \hbox{if $\psi(e) < 0$.}
  \end{array}
\right.
\]
Note that $\psi'$ preserves the signs of charges assigned by $\psi$, it merely `scales down' the positive charges. The following gives an important connection between $(T_{a,b,c})$-covered assignments and triangle-covered ones.

\begin{lemma}\label{L:reductionTriCov}
Let $a,b,c \in \mathbb{R}^+$ such that $a \geq b \geq c$, let $\alpha \in (0,1]$, and let $\psi$ be a $(T_{a,b,c})$-covered charge assignment to a graph $G$. If $\alpha \geq \max\{\frac{2c}{a+b},\frac{b+c}{2a}\}$, then the $\alpha$-reduction of $\psi$ is triangle-covered.
\end{lemma}

\begin{proof}
Let $\psi'$ be the $\alpha$-reduction of $\psi$. Let $e_1$, $e_2$, and $e_3$ 
be three distinct edges of $G$ that form a triangle, and let $\sigma' \coloneqq \psi'(e_1)+\psi'(e_2)+\psi'(e_3)$. We must show that $\sigma' \geq 0$. Suppose without loss of generality that $\psi(e_1) \leq \psi(e_2) \leq \psi(e_3)$. Because $\psi$ is $(T_{a,b,c})$-covered, we have
\begin{equation}\label{E:psiCovered}
    a\psi(e_1)+b\psi(e_2)+c\psi(e_3) \geq 0.
\end{equation}

\noindent \textbf{Case 1: no negative edge.}  Suppose $\psi(e_1) \geq 0$. Then $\sigma' \geq 0$ because all three terms are nonnegative. \smallskip

\noindent \textbf{Case 2: one negative edge.} Suppose $\psi(e_1) < 0$ and $\psi(e_2) \geq 0$. Using first $b \geq c$ and $0 \leq \psi(e_2) \leq \psi(e_3)$, and then $\alpha \geq \frac{b+c}{2a}$, we obtain
\[0 \stackrel{\cref{E:psiCovered}}{\leq} a\psi(e_1)+\tfrac{b+c}{2}\bigl(\psi(e_2)+\psi(e_3)\bigr) \leq  a\psi'(e_1)+\tfrac{b+c}{2}\cdot\tfrac{1}{\alpha}\bigl(\psi'(e_2)+\psi'(e_3)\bigr) \leq a\sigma'\]
and hence $\sigma' \geq 0$. \smallskip

\noindent\textbf{Case 3: two negative edges.} Suppose $\psi(e_2) < 0$. Then $\psi(e_3) \geq 0$ by \cref{E:psiCovered}. Hence, using first $a \geq b$ and $\psi(e_1) \leq \psi(e_2) \leq 0$, and then $\alpha \geq \frac{2c}{a+b}$, we obtain
\[0 \stackrel{\cref{E:psiCovered}}{\leq} \tfrac{a+b}{2}\bigl(\psi(e_1)+\psi(e_2)\bigr)+c\psi(e_3) =  \tfrac{a+b}{2}(\psi'(e_1)+\psi'(e_2))+c\cdot\tfrac{1}{\alpha}\psi'(e_3) \leq \tfrac{a+b}{2}\sigma'.\]
So again $\sigma' \geq 0$.\smallskip
\end{proof}

Finally we need the observation made in \cite{DelPos2026} that we can make small perturbations to any charge assignment in such a way that none of its important properties are altered but we remove any zero charges or ties. A charge assignment $\phi$ to a graph $G$ is \emph{nonzero} if $\phi(e) \neq 0$ for all $e \in E(G)$ and is \emph{untied} if $\phi(e) \neq \phi(f)$ for all distinct $e,f \in E(G)$.

\begin{lemma}\label{L:nonzeroUntied}
Let $a,b,c \in \mathbb{R}^+$. If there is a $(T_{a,b,c})$-covered charge assignment to a graph $G$ with negative total weight, then there is a nonzero, untied and $(T_{a,b,c})$-covered charge assignment to $G$ with negative total weight.
\end{lemma}

\begin{proof}
Let $\phi$ be a $(T_{a,b,c})$-covered charge assignment to $G$ with total weight $-t$ where $t>0$. Let $C\coloneqq \{\phi(e):e \in E(G)\} \cup \{0\}$ and let $\varepsilon\coloneqq \min\{|c-c'|:c,c' \in C, c \neq c'\}$ (this is well defined since $C$ must contain 0 and a negative element). Let $m \coloneqq |E(G)|$, let $(e_1,\ldots,e_m)$ be an arbitrary ordering the edges of $G$, and choose distinct elements $\epsilon_1,\ldots,\epsilon_m$ of the real interval $(0,\min\{\epsilon,\frac{t}{m}\})$. Define a charge assignment $\phi'$ to $G$ by $\phi'(e_i)\coloneqq \phi(e_i)+\epsilon_i$ for each $i \in \{1,\ldots,m\}$. Then $\phi'$ is $(T_{a,b,c})$-covered because $\phi'(e) > \phi(e)$ for all $e \in E(G)$, $\phi'$ is nonzero and untied because $\epsilon_1,\ldots,\epsilon_m$ are distinct and each is
in $(0,\varepsilon)$, and $\phi'$ has negative total weight because $\sum_{i=1}^m \varepsilon_i < m \cdot \frac{t}{m} = t$.
\end{proof}

\section{Discharging} \label{sec:Discharging}

We now introduce a slight generalisation of the discharging method used in \cite{DelPos2026}. Let $G$ be a graph and $\phi$ be a nonzero untied charge assignment to $G$. We denote the sets of positive and negative edges by $E^+$ and $E^-$ respectively. For an edge $e$ of $G$, we define $T^+(e)$ and $T^-(e)$ to be the set of all edges of $G$ that are in a triangle with $e$ and are, respectively, positive and negative. These notations are implicitly dependent on the graph $G$ and charge assignment $\phi$ but we use them unambiguously, sometimes taking advantage of the fact that sets of positive and negative edges under a given charge assignment and a reduction of it are the same. Let $f \in E^-$, let $e \in T^+(f)$, and let $d$ be the edge that forms a triangle with $e$ and $f$. Let $\tau \in \mathbb{R}^+$ be a given \emph{sponsorship threshold}. When $\phi(d) >0$, we say $e$ is a \emph{$\tau$-sponsor} of $f$ if $\phi(e)>\phi(d)$ and $\tau \phi(e) \geq |\phi(f)|$. We define the \textit{demand of $f$ on $e$} by
\[
\dem^{\tau}_{\phi}(e,f) \coloneqq
\begin{cases} 
      |\phi(f)| & \text{if $\phi(d)<0$} \\
      |\phi(f)| & \text{if $e$ is a $\tau$-sponsor of $f$} \\
      0 & \text{if $d$ is a $\tau$-sponsor of $f$} \\
      \frac{\phi(e)}{\phi(e)+\phi(d)}|\phi(f)| & \text{otherwise.}
\end{cases}
\]
Further, given $f \in E^-$ and $e \in T^+(f)$, we define the \textit{flow of charge from $e$ to $f$} to be
\[\ch^{\tau}_{\phi}(e,f) \coloneqq
\begin{cases}
    0 & \text{if $\dem^{\tau}_{\phi}(e,f)=0$}\\
    \frac{\dem^{\tau}_{\phi}(e,f)}{\sum_{f' \in T^-(e)}\dem^{\tau}_{\phi}(e,f')} \cdot \phi(e) & \text{otherwise.}
\end{cases}\]
Finally, define the \textit{final charge on an edge $d \in E(G)$} to be
\[\fch^{\tau}_{\phi}(d) \coloneqq
\begin{cases} 
      \phi(d) - \sum_{f \in T^-(d)}\ch^{\tau}_{\phi}(d,f) & \text{if $\phi(d)>0$} \\
      \phi(d) + \sum_{e \in T^+(d)}\ch^{\tau}_{\phi}(e,d) & \text{if $\phi(d)<0$.}
\end{cases}
\]

Observe that this means that $\fch^{\tau}_{\phi}(e) \geq 0$ for every $e \in E(G)$ with $\phi(e)>0$ and the sum of the final charges $\sum_{d\in E(G)}\fch^{\tau}_{\phi}(d)$ is equal to the (original) total charge of $\phi$. So, to show that the total charge of $\phi$ is nonnegative, it suffices to show that there exists some $\tau \in \mathbb{R}^+$ such that $\fch^{\tau}_{\phi}(f)$ is positive for all $f \in E^-$. In all this notation, we omit $\tau$ when it is equal to 1, as is always the case in \cite{DelPos2026}.
The following facts follow from these rules.

\begin{lemma}\label{L:scaled}
Let $\alpha \in (0,1]$, let $G$ be a graph, let $\psi$ be a nonzero and untied charge assignment to $G$, and let $\psi'$ be the $\alpha$-reduction of $\psi$. For all $f\in E^-$, the following hold.
\begin{enumerate}
    \item
$\dem_{\psi'}(e,f)=\dem^{\alpha}_{\psi}(e,f)$ for all $e \in T^+(f)$.
    \item
$\ch_{\psi'}(e,f)=\alpha\,\ch^{\alpha}_{\psi}(e,f)$ for all $e \in T^+(f)$.
    \item
If $\fch_{\psi'}(f)\geq \psi'(f)+\alpha|\psi'(f)|$, then $\fch^{\alpha}_{\psi}(f) \geq 0$.
\end{enumerate}
\end{lemma}

\begin{proof}
Let $f \in E^-$ and let $e \in T^+(f)$. Part (i) follows from the definition of $\dem^\alpha_\phi(e,f)$ because $e$ is an $\alpha$-sponsor of $f$ under $\psi$ if and only if $e$ is a $1$-sponsor of $f$ under $\psi'$. Part (ii) then follows from part (i), the definition of $\ch^{\alpha}_{\phi}(e,f)$, and the fact that $\psi'(e)=\alpha\,\psi(e)$. In turn, part (iii) follows from part (ii), the definition of $\fch^{\alpha}_{\phi}(f)$, and the fact that $\psi'(f)=\psi(f)$.
\end{proof}

The key result from \cite{DelPos2026} that we wish to use is the following.

\begin{lemma}[{\cite[Lemma 10.4]{DelPos2026}}]\label{L:DelPosPacking} Let $\gamma \in [0,\frac{1}{4}]$. Let $\phi$ be a nonzero, untied and triangle-covered charge assignment to an $n$-vertex graph $G$ with minimum degree at least $(\frac{3}{4}-\gamma)n$. If $f$ is a negative edge of $G$, then $\fch_{\phi}(f) \geq -4\gamma|\phi(f)|$.
\end{lemma}

We are now in a position to prove our \cref{T:weightedTriDecomp}. Essentially, \cref{L:reductionTriCov} and \cref{L:DelPosPacking} can be combined to give the upper bound on $\delta^*_{T_{a,b,c}}$ given in part (i). Part (ii) is merely noting when this upper bound matches the lower bound given by \cref{C:LBforTris}.

\begin{proof}[\textup{\textbf{Proof of \cref{T:weightedTriDecomp}.}}]
We first prove (i). Suppose otherwise. By \cref{L:FarkasForConds} there must be a $(T_{a,b,c})$-covered charge assignment $\psi$ to $G$ with negative total charge. By \cref{L:nonzeroUntied}, we may assume that $\psi$ is nonzero and untied. Let $\alpha \coloneqq \max\{\frac{2c}{a+b},\frac{b+c}{2a}\}$. We will find a contradiction by showing $\fch^{\alpha}_{\psi}(f) \geq 0$ for all $f \in E^-$.

Let $\psi'$ be the $\alpha$-reduction of $\psi$. Then $\psi'$ is nonzero and untied by definition and is triangle-covered by \cref{L:reductionTriCov}. Let $\gamma \coloneqq \frac{1-\alpha}{4}$ and note that $\frac{3}{4}-\gamma=\frac{1}{2}+\frac{\alpha}{4}=\frac{1}{2}+\max\{\frac{c}{2(a+b)},\frac{b+c}{8a}\}$, so $G$ has minimum degree at least $(\frac{3}{4}-\gamma)n$ by hypothesis. Let $f \in E^-$. Applying \cref{L:DelPosPacking} to $\psi'$ with our chosen value of $\gamma$, we see that 
\[\fch_{\psi'}(f) \geq -4\gamma|\psi'(f)|=\psi'(f)+(1-4\gamma)|\psi'(f)| =  \psi'(f)+\alpha|\psi'(f)|.\]
Thus, $\fch^{\alpha}_{\psi}(f) \geq 0$ by \cref{L:scaled}(iii) and the proof of (i) is complete. Noting that $c \geq \frac{b(a+b)}{3a-b}$ is equivalent to $\frac{c}{2(a+b)} \geq \frac{b+c}{8a}$, (ii) follows directly from (i) and \cref{C:LBforTris}.  
\end{proof}

Our remaining results now follow quite easily.

\begin{proof}[\textup{\textbf{Proof of \cref{T:approx}}}]
\cref{T:weightedTriDecomp}(i) implies $\delta^*_{T_{a,b,c}} \leq \frac{1}{2}+\max\{\frac{c}{2(a+b)},\frac{b+c}{8a}\}$. Thus (i) follows from \cref{T:condensation} since $T_{a,b,c}$ is a condensation of $F$. Now suppose $F$ obeys the hypotheses of (ii). Since $c \geq \frac{b(a+b)}{3a-b}$, we have $\frac{c}{2(a+b)} \geq \frac{b+c}{8a}$ and hence $\delta^{0+}_F \leq \frac{1}{2}+\frac{c}{2(a+b)}$. It remains to show that $\delta^{0+}_F \geq \frac{1}{2}+\frac{c}{2(a+b)}$. Let $\epsilon>0$ be arbitrarily small such that $\delta \coloneq \frac{1}{2}+\frac{c}{2(a+b)}-\epsilon$ is rational, and let $n$ be a positive even integer such that $\delta n$ is also an integer. Consider a graph $G$ in $\mathcal{B}_n(\delta)$ (defined after \cref{C:LBforTris}) and let $G'$ be a subgraph of $G$ that has an $F$-decomposition. We will show that $|E(G)|-|E(G')| \geq \epsilon n^2$ to establish $\delta_F^{0+} \geq \frac{1}{2}+\frac{c}{2(a+b)}$ and so complete the proof.

Note that $G$ has $m_0 \coloneq (\frac{\delta}{2}-\frac{1}{4})n^2$ edges not in its $K_{\frac{n}{2},\frac{n}{2}}$ subgraph. Because no graph formed by removing fewer than $c$ edges from $F$ is bipartite, each copy of $F$ in $G$ uses at least $c$ of these edges. Hence, $|E(G')| \leq \frac{a+b+c}{c}\,m_0=(1+\frac{a+b}{c})m_0$. So
\[\mfrac{|E(G)|-|E(G')|}{n^2} \geq \mfrac{\delta}{2} - \left(\mfrac{a+b}{c}+1\right)\left(\mfrac{\delta}{2}-\mfrac{1}{4}\right) = \mfrac{1}{4}-\mfrac{(a+b)(2\delta-1)}{4c}=\mfrac{a+b}{2c}\epsilon \geq \epsilon. \qedhere\]
\end{proof}

\begin{proof}[\textup{\textbf{Proof of \cref{C:multi}}}]
We can apply \cref{T:approx} with $F\coloneq K_{x,y,z}$, $a \coloneq xy$, $b \coloneq xz$ and $c \coloneq yz$. A routine calculation shows that $c \geq \frac{b(a+b)}{3a-b}$ is equivalent to $z \geq \frac{y(3y-z)}{y+z}$. Further, noting that there are $xyz=xc$ triangles in $K_{x,y,z}$, and deleting an edge destroys at most $x$ of these, we can confirm the hypotheses of \cref{T:approx}(ii) hold and hence the result follows.
\end{proof}

\begin{proof}[\textup{\textbf{Proof of \cref{C:wheel}}}]
Consider a 3-colouring of $W_n$ in which one colour class contains only the dominating vertex and the other two colour classes each contain $\frac{n-1}{2}$ of the other vertices. This shows we can apply \cref{T:approx} with $F \coloneq W_n$, $a \coloneq n-1$, $b\coloneq\frac{n-1}{2}$ and $c\coloneq\frac{n-1}{2}$. Further, $W_n$ contains $n-1=2c$ triangles and observe that removing an edge destroys at most two of these. So the hypotheses of \cref{T:approx}(ii) hold and the result follows.
\end{proof}

\begin{proof}[\textup{\textbf{Proof of \cref{T:main}.}}]
Consider a 3-colouring of $C_\ell$ in which one colour class contains exactly one vertex and the other two colour classes each contain $\frac{\ell-1}{2}$ of the other vertices. Then it is easy to check that we can apply \cref{T:approx}(ii) with $F\coloneq C_\ell$, $a \coloneq \ell-2$, $b \coloneq 1$ and $c \coloneq 1$ to obtain $\delta^{0+}_{C_{\ell}} = \frac{1}{2}+\frac{1}{2(\ell-1)}$. Now the result follows by \cref{T:approxToIntegral}.
\end{proof}

\section{Limitations of the discharging rules}\label{S:limitations}

We observed in \cref{T:weightedTriDecomp}(ii) that the general upper bound $\delta_{T_{a,b,c}}^*\leq \frac{1}{2}+\max\{\frac{c}{2(a+b)},\frac{b+c}{8a}\}$ obtained in \cref{T:weightedTriDecomp}(i) matches the lower bound $\delta_{T_{a,b,c}}^*\geq \max\{B(a,b,c),U(a,b,c)\}$ from  \cref{C:LBforTris} whenever $a\geq b\geq c \geq \frac{b(a+b)}{3a-b}$. When $c<\frac{b(a+b)}{3a-b}$, it is not too difficult to verify that $\frac{1}{2}+\frac{b+c}{8a}> \max\{B(a,b,c),U(a,b,c)\}$, and so we do not determine $\delta_{T_{a,b,c}}^*$ for this family of weighted triangles.

However, this restriction on $(a,b,c)$ is not merely an artefact of our analysis; \cref{T:weightedTriDecomp}(i) cannot be improved without developing new discharging rules, or modifying the proof strategy altogether. Recall that the general idea is that, by \cref{L:FarkasForConds}, it suffices to show that any $(T_{a,b,c})$-covered charge assignment has nonnegative total charge. To do so, we use the discharging rules of \cite{DelPos2026} to show that, given any $(T_{a,b,c})$-covered charge assignment $\phi$, the charges can be redistributed to obtain another $(T_{a,b,c})$-covered charge assignment $\phi'$, with the same total charge, but where all charges are nonnegative.

In the following lemma, however, we show for  $\delta$ below but arbitrarily close to $\frac{1}{2}+\frac{b+c}{8a}$, that there exists a graph of minimum degree $\delta n$ and a $(T_{a,b,c})$-covered charge assignment $\phi$ such that, after applying the discharging rules, there is still an edge of negative charge.
This demonstrates the necessity of introducing new rules, or of developing a different proof strategy which enables charges to be counted globally rather than individually. The graphs we employ are in the family $\mathcal{U}_n$ (defined after \cref{C:LBforTris}).

\begin{lemma}\label{L:fail}
Let $a,b,c \in \mathbb{R}^+$ such that $a \geq b \geq c$ and let $\delta \in (\frac{1}{2}, \frac{1}{2}+\frac{b+c}{8a})$ be rational. For all $n_0 \in \mathbb{N}$, there is $n\geq n_0$, an $n$-vertex graph $G$ with $\delta(G)=\delta n$, and a nonzero, untied and $(T_{a,b,c})$-covered charge assignment $\phi$ to $G$ such that, for all $\tau \in \mathbb{R}^+$, we have $\fch^\tau_\phi(f) < 0$ for some $f \in E^-$.    
\end{lemma}

\begin{proof}
Let $n$ be a positive even integer such that $\delta n$ is also an integer. Let $G$ be the complete $4$-partite graph with parts $A_1,A_2,B_1,B_2$ where $|A_1|=|A_2|=(\delta-\frac{1}{2})n$ and $|B_1|=|B_2|=(1-\delta)n$. Let $A\coloneqq A_1 \cup A_2$ and $B\coloneqq B_1 \cup B_2$. Let $\phi^*$ be a charge assignment of $1$ to each edge with one end in $A$ and the other in $B$, and of $-\frac{b+c}{a}$ to each other edge. Now we untie $\phi^*$. Let $m \coloneqq |E(G)|$. Since $\delta<\frac{1}{2}+\frac{b+c}{8a}$, we have $4(2\delta-1)<\frac{b+c}{a}$ and hence we can choose $\epsilon \in \mathbb{R}^+$ sufficiently small such that $\epsilon+4(2\delta-1)(1+\epsilon)<\frac{b+c}{a}$. Choose a strictly increasing sequence $(\epsilon_1,\ldots,\epsilon_{m})$ of positive real numbers such that $\epsilon_{m}<\epsilon$. Choose an ordering $(e_1,\ldots,e_m)$ of the edges of $G$ uniformly at random. For each $i \in \{1,\ldots,m\}$, let $\phi(e_i)\coloneqq \phi^*(e_i)+\epsilon_i$. Note that $\phi^*$ is $(T_{a,b,c})$-covered and hence so is $\phi$. Fix $\tau \in  \mathbb{R}^+$ and an edge $f_0$ with both endpoints in $B$. To complete the proof, we show that $\mathbb{E}\,\fch^\tau_\phi(f_0) < 0$ and hence that $\fch^\tau_\phi(f_0) < 0$ for some choice of $(e_1,\ldots,e_m)$.

Let $e \in E^+$, let $f,f' \in T^-(e)$, and let $d$ and $d'$ be the edges that form triangles with $e$ and $f$ and with $e$ and $f'$ respectively. Note that by construction $e$, $d$ and $d'$ each have one endpoint in $A$ and the other in $B$, while each of $f$ and $f'$ has both its endpoints in $A$ or both in $B$. In particular, $\phi^*(f)=\phi^*(f')=-\frac{b+c}{a}$ and $\phi^*(d)=\phi^*(d')=1$. Since the ordering $(e_1, \dots, e_m)$ is chosen uniformly at random, the probability distributions of $(\phi(f),\phi(d))$ and $(\phi(f'),\phi(d'))$ are identical. Therefore,
\[\mathbb{E}\, \mfrac{\dem^\tau_\phi(e,f')}{\dem^\tau_\phi(e,f)}=1,\]
and since $|T^-(e)|=\frac{n}{2}$, we obtain 
\[\mathbb{E}\,\ch^\tau_\phi(e,f) < \mathbb{E}\,\mfrac{\dem^{\tau}_{\phi}(e,f)}{\sum_{f' \in T^-(e)}\dem^{\tau}_{\phi}(e,f')} \cdot (1+\epsilon) = \mathbb{E}\,\left(\sum_{f' \in T^-(e)}\mfrac{\dem^{\tau}_{\phi}(e,f')}{\dem^{\tau}_{\phi}(e,f)}\right)^{-1} \cdot (1+\epsilon) = \mfrac{2}{n}(1+\epsilon).\]
Hence, as $|T^+(f_0)|=2(2\delta-1)n$, we have
\[\mathbb{E}\, \fch^\tau_\phi(f_0)<\phi(f_0)+4(2\delta-1)(1+\epsilon)<-\mfrac{b+c}{a}+\epsilon+4(2\delta-1)(1+\epsilon) < 0,\]
where the final inequality follows by our choice of $\epsilon$. 
\end{proof}

\section{Conclusion}

In this note, we showed how the `discharging in the dual' and `condensation' methods can be used to bound the fractional decomposition threshold of weighted triangles.
Weighted triangles $T_{a,b,c}$ can be visualised in \cref{fig:feasible_mu} below.
Assuming $a \geq b \geq c$, the full range of possible values of $(a,b,c)$ are represented in the shaded areas of the plot. The dark grey region represents the triangles for which \cref{T:weightedTriDecomp}(ii) fully determines the fractional decomposition threshold $\delta_{T_{a,b,c}}^*$.
In the light grey region however, the upper bound on $\delta_{T_{a,b,c}}^*$ given by \cref{T:weightedTriDecomp}(i) does not match the lower bound given by \cref{C:LBforTris}. In that region, there are two types of triangles. First, there are triangles where $B(a,b,c) \geq U(a,b,c)$ but \cref{T:weightedTriDecomp}(ii) does not apply; examples include $T_{4,4,3}$ and $T_{4,2,1}$. Second, there are triangles where $B(a,b,c) < U(a,b,c)$ (note \cref{T:weightedTriDecomp}(ii) cannot apply to any of these); examples include $T_{2,2,1}$ and $T_{4,3,1}$.

It is interesting to note that, when $B(a,b,c) > U(a,b,c)$ but \cref{T:weightedTriDecomp}(ii) does not apply, the better lower bound on $\delta^*_{T_{a,b,c}}$ is proved with graphs in the family $\mathcal{B}_n$ (see the proof of \cref{C:LBforTris}), while the proof that the discharging methods of \cite{DelPos2026} fail uses a graph in the family $\mathcal{U}_n$ (see the proof of \cref{L:fail}).

\begin{figure}[h]
    \centering
\includegraphics{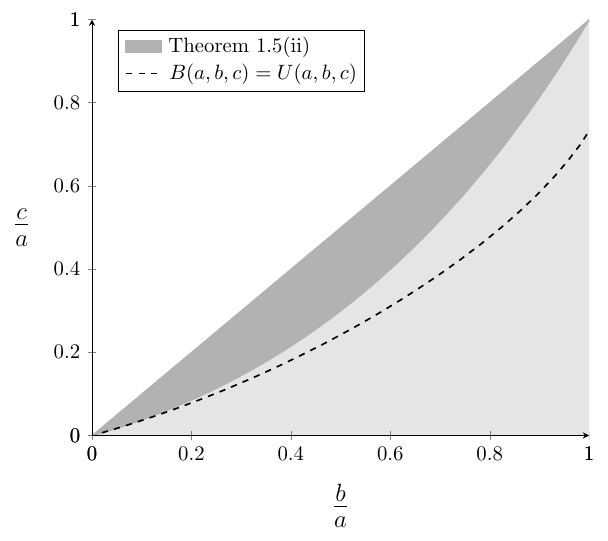}
\caption{Representation of the possible weighted triangles $T_{a,b,c}$, with $a\geq b\geq c$. In the dark grey
region, thresholds are fully determined by \cref{T:weightedTriDecomp}(ii). In the light grey region thresholds are not determined. Above the dashed line, $B(a, b, c) > U (a, b, c)$. }
    \label{fig:feasible_mu}
\end{figure}

As \cref{T:main,T:approx} indicate, finding fractional decomposition thresholds for weighted triangles can lead to obtaining decomposition thresholds, or approximate decomposition thresholds, for 3-chromatic graphs. So it is natural to wish to generalise \cref{T:weightedTriDecomp} to a wider range of weighted triangles. We showed in \cref{S:limitations}, however, that to do so using discharging in the dual would require developing different discharging rules to those in \cite{DelPos2026}. 

In the case of odd cycles, we could derive the decomposition threshold from the approximate one due to the absorption results of \cite{BarKuhLoOst}. For other $3$-chromatic graphs such absorption results are not yet available. A general upper bound is given in \cite{GloKuhLoMonOst}.

\section*{Acknowledgments}

Thanks to Michelle Delcourt and Luke Postle for valuable discussions.

\section*{AI use statement}

No AI was used in the preparation of this manuscript.

\bibliography{refs}
\bibliographystyle{abbrv}

\end{document}